# Validity of Borodin and Kostochka Conjecture for $4K_1$–free Graphs

Medha Dhurandhar

**Abstract:** Problem of finding an optimal upper bound for the chromatic no. of even $3K_1$-free graphs is still open and pretty hard. Here we prove **Borodin & Kostochka Conjecture** for $4K_1$-free graphs G i.e. if $\Delta(G) \geq 9$ and G is $4K_1$-free, then $\chi(G) \leq \max\{\omega, \Delta-1\}$.

**Introduction:**

In [1], [2], [3], [4] chromatic bounds for graphs are considered especially in relation with $\omega$ and $\Delta$. Gyárfás [5] and Kim [6] show that the optimal $\chi$-binding function for the class of $4K_1$-free graphs has order $\omega^2/\log(\omega)$. If we forbid additional induced subgraphs, the order of the optimal $\chi$-binding function drops below $\omega^2/\log(\omega)$. In 1941, Brooks' theorem stated that for any connected undirected graph $G$ with maximum degree $\Delta$, the chromatic number of $G$ is at most $\Delta$ unless $G$ is a complete graph or an odd cycle, in which case the chromatic number is $\Delta + 1$ [5]. In 1977, **Borodin & Kostochka** [6] conjectured that if $\Delta(G) \geq 9$, then $\chi(G) \leq \max\{\omega, \Delta-1\}$. In 1999, Reed proved the conjecture for $\Delta \geq 10^{14}$ [7]. Also D. W. Cranston and L. Rabern [8] proved it for claw-free graphs. Here we prove **Borodin & Kostochka** conjecture for $4K_1$-free graphs.

**Notation:** For a graph G, V(G), E(G), $\Delta$, $\omega$, $\chi$ denote the vertex set, edge set, maximum degree, size of a maximum clique, chromatic number of G resply. For $u \in V(G)$, $N(u) = \{v \in V(G) / uv \in E(G)\}$, and $\overline{N(u)} = N(u) \cup (u)$. If $S \subseteq V$, then $<S>$ denotes subgraph of G induced by S. If C is some coloring of G and if $u \in V(G)$ is colored m in C, then u is called a m-vertex, if N(u) has a unique r-vertex, then r is called a unique color of u and if N(u) has more than one r-vertex, then r is called a repeat color of u. Also if P is a path in G s.t. vertices on P are alternately colored say i and j, then P is called an i-j path. All graphs considered henceforth are simple. We consider here simple and undirected graphs. For terms which are not defined herein we refer to Bondy and Murty [9].

**Main Result:** Let G be $4K_1$-free and $\Delta \geq 9$, then $\chi \leq \max\{\Delta-1, \omega\}$.

Proof: Let if possible G be a smallest, connected, $4K_1$-free graph with $\Delta \geq 9$ and $\chi > \max\{\Delta-1, \omega\}$. Then clearly as $G \neq C_{2n+1}$ or $K_{|V(G)|}$, $\chi = \Delta > \omega$. Let $u \in V(G)$. Then $G-u \neq K_{|V(G)|-1}$ (else $\chi = \omega$). **If $\Delta(G-u) \geq 9$**, then by minimality $\chi(G-u) \leq \max\{\omega(G-u), \Delta(G-u)-1\}$. Clearly if $\omega(G-u) \leq \Delta(G-u)-1$, then $\chi(G-u) = \Delta(G-u)-1 \leq \Delta-1$ and otherwise $\chi(G-u) = \omega(G-u) \leq \omega < \Delta$. In any case $\chi(G-u) \leq \Delta-1$. **Also if $\Delta(G-u) < 9$**, then as $G-u \neq C_{2n+1}$, by Brook's Theorem $\chi(G-u) \leq \Delta(G-u) < 9 \leq \Delta$. Thus always $\chi(G-u) \leq \Delta-1$ and in fact, $\chi(G-u) = \Delta-1$ and deg $v \geq \Delta-1$ $\forall v \in V(G)$.

Let $u \in V(G)$ be s.t. deg $u = \Delta$. Let $S = \{1,..., \Delta-1, \Delta\}$ be a $\Delta$-coloring of G with only u colored $\Delta$. Then N(u) has $\Delta-2$ vertices $A_i$ with unique colors i ($1 \leq i \leq \Delta-2$) and a pair of vertices say X, Y with the same color $\Delta-1$. Clearly $A_i$ has a j-vertex for $1 \leq i \neq j \leq \Delta-2$ (else color $A_i$ by j, u by i).

**Case 1:** $\exists$ a $(\Delta-1)$-coloring of G-u s.t. $A_iA_j \notin E(G)$ for some i, j $\in \{1,..., \Delta-2\}$.

**(A)** For no m, $A_m$ is the only m-vertex of both $A_i$ and $A_j$ for $1 \leq i, j, m \leq \Delta-2$.
Let if possible $A_i$, $A_k$ both have $A_m$ as the only m-vertex. Then as $A_m$ has at the most one repeat color, w.l.g. $A_j$ be the only j-vertex of $A_m$. Then color $A_i$, $A_j$ by m, $A_m$ by j, u by i, a contradiction.
**(B)** $A_i$, $A_j$ do not have more than two common adjacent $A_k$s in N(u).
Let $A_i$, $A_j$ be both adjacent to say $A_k$, $A_l$, $A_m$ $1 \leq i, j, k, l, m \leq \Delta-2$. As each of $A_i$, $A_j$ has at the most one repeat color, w.l.g. let $A_m$ be the only m-vertex of both $A_i$ and $A_j$, a contradiction to **(A)**.
**(C)** $A_i$ is non-adjacent to at the most three $A_k$ s $\Rightarrow$ As $\Delta-2 \geq 7$, $A_i$ is adjacent to at least three $A_m$, $1 \leq i, k, m \leq \Delta-2$.
Let if possible $A_1A_k \notin E(G)$ for $2 \leq k \leq 5$. As G is $4K_1$-free, $\exists$ at most two more 1-vertices $a_{11}$, $a_{12}$ and as $\exists$ 1-k path from $A_1$ to $A_k$, either $a_{11}$ or $a_{12}$ is adjacent to $A_k$ with two k-vertices for $2 \leq k \leq 5$. Again $a_{1i}$ cannot have three repeat colors (else $N(a_{1i})$ has a color say r missing. Color $a_{1i}$ by r. Then either **(i)** some $A_k$ ($2 \leq k \leq 5$) has no 1-vertex, hence color $A_k$ by 1, u by k or **(ii)** $a_{12}A_k \in E(G)$ ($2 \leq k \leq 5$), $a_{12}$ has four repeat colors and $N(a_{12})$ has color t missing. Color $a_{12}$ by t, $A_k$ by 1, u by k). Thus w.l.g. let $A_ia_{11}$, $A_ja_{12} \in E(G)$ for i = 2, 3 and j = 4, 5 s.t. $a_{11}$ has two repeat colors 2, 3 and $a_{12}$ has two repeat colors 4, 5. Clearly $A_1A_j \in E(G)$ $\forall$ $6 \leq i \leq \Delta-2$ (else either $a_{11}$ or $a_{12}$ has three repeat colors).

**Claim 1:** Whenever $A_1$ has a unique i-vertex say B for $6 \leq i \leq \Delta-1$, $A_1$ is the only 1-vertex of B.
Let if possible $Ba_{11} \in E(G)$. Then B has a unique m-vertex for $2 \leq m \leq \Delta-1$ (else N(B) has some color r missing. Color B by r, $A_1$ by i, u by 1). As $a_{11}$ has two repeat colors 2, 3, B is its only i-vertex. Then G has at the most one more i-vertex say b (else $<A_1, a_{11}, b_{11}, b_{12}> = 4K_1$). Again by **(A)**, B is not the only i–vertex of any $A_k$, for $2 \leq k \leq 5$. Hence $A_kb \in E(G)$ for $2 \leq k \leq 5$. Now $A_kB \notin E(G)$ for k = 2, 3 (else color $A_k$ by 1, $a_{11}$ by i, B by 2/3, $A_1$ by i, u by 1) and b has two k-vertices for k = 2, 3 (else color $A_k$ by i, b by k, u by k) $\Rightarrow A_m$ is the only m-vertex of b for m = 4, 5 (else b has color r missing in N(b). Color b by r, $A_2$ by i, u by 2). Now $A_m$ has two i-vertices (else b is the only i-vertex of $A_m$. Color b by m, $A_m$ by i, u by m), m $\in$ {4, 5} $\Rightarrow a_{12}$ is the only 1-vertex of $A_m$, m $\in$ {4, 5}. Again $Ba_{12} \notin E(G)$ (else B has three 1-vertices and color say r missing in N(B). Color B by r, $A_1$ by i, u by 1) $\Rightarrow a_{12}b \in E(G)$. Then color b by 4, $A_4$ by 1, $a_{12}$ by i, u by 4, a contradiction. This proves **Claim 1**.

Now $a_{1k}$ has an i-vertex for k = 1, 2 (else color $a_{1k}$ by i. If $a_{1k}$ is the only 1-vertex of $A_m$ ($2 \leq m \leq 5$), then color $A_m$ by 1, u by m and if every $A_m$ has two 1-vertices, then if k = 1 (2), color $A_2$ ($A_4$) by 1, $a_{12}$ ($a_{11}$) by 2 (4), u by 2 (4)).

Let $a_{11}b_{i1} \in E(G)$. As $a_{11}$ has two repeat colors 2, 3, $b_{i1}$ is the only i-vertex of $a_{11}$.

**Claim 2:** $a_{11}$ is the only 1-vertex of $b_{i1}$.
Let if possible $a_{12}b_{i1} \in E(G)$. As $a_{12}$ has two repeat colors 4, 5, $b_{i1}$ is the only i-vertex of $a_{12}$. Then G doesn't have an i-vertex say $b_{12} \notin \{B, b_{11}\}$ (else $<a_{11}, a_{12}, B, b_{12}> = 4K_1$). Again by **(A)**, B cannot be the only i-vertex of any $A_m$ for $2 \leq m \leq 5$. Hence $A_mb_{i1} \in E(G)$ for $2 \leq m \leq 5$. **If $A_k$ is the only k-vertex of $b_{i1}$** for some k, $2 \leq k \leq 5$, then if $b_{i1}$ is the only i-vertex of $A_k$, color $A_k$ by 1, $b_{i1}$ by k, u by k and if $A_k$ has two i-vertices, then $a_{1j}$ being the only 1-vertex of $A_k$, color $A_k$ by 1, $a_{1j}$ by i, $b_{i1}$ by k, u by k, contradictions in both the cases. **Hence let $b_{i1}$ have repeat colors k $\forall$ k, $2 \leq k \leq 5$**. But then $b_{i1}$ has color r missing in $N(b_{i1})$. Color $b_{i1}$ by r and $a_{11}$ by i. Then $A_2a_{12} \in E(G)$ (else color $A_2$ by 1, u by 2). Again as $a_{12}$ has two repeat colors 4, 5, $A_2$ is its only 2-vertex and hence color $A_2$ by 1, $a_{12}$ by 2, u by 2, a contradiction. This proves **Claim 2**.

Similarly if $b_{i2}$ is an i-vertex of $a_{12}$, then $a_{12}$ ($b_{i2}$) is the only 1-vertex (i-vertex) of $b_{i2}$ ($a_{12}$). Now $A_mb_{i1}$, $A_nb_{i2} \in E(G)$ for m = 2, 3 and n = 4, 5 (else let $A_2b_{i1} \notin E(G)$. If $a_{11}$ is the only 1-vertex of $A_2$, then color $a_{11}$ by i, $b_{i1}$ by 1, $A_2$ by 1, u by 2 and if $A_2a_{12} \in E(G)$, then color $a_{11}$ by i, $b_{i1}$ by 1, $a_{12}$ by 2, $A_2$ by 1, u by 2).

As $A_1$ has at the most one repeat color, w.l.g. let $A_1$ have unique 2, 3, 4 vertices. Let P (R) be a 2-1 (4-1) path from $A_2$ ($A_4$) to $A_1$. As $a_{12}$ ($a_{11}$) has a unique 2-vertex (4-vertex), clearly P = $\{A_2, a_{11}, a_{21}, A_1\}$ and R = $\{A_4, a_{12}, a_{41}, A_1\}$.

**Claim 3:** $a_{21}a_{12}$, $A_2a_{12} \notin E(G)$. Similarly $a_{41}a_{11}$, $A_4a_{11} \notin E(G)$.
Let if possible $a_{21}a_{12} \in E(G)$. Then G has no other 2-vertex $a_{22} \notin \{A_2, a_{21}\}$ (else $<a_{22}, a_{12}, A_1, A_2> = 4K_1$). Also $a_{21}b_{i1} \in E(G)$ (else $A_2$ is the only 2-vertex of $b_{i1}$. If $b_{i1}$ is the only i-vertex of $A_2$, then color $b_{i1}$ by 2, $A_2$ by i, u by 2 and if $A_2$ has two i-vertices, then color $b_{i1}$ by 2, $A_2$ by 1, $a_{11}$ by i, u by 2). As $a_{21}$ has three 1-vertices, $Ba_{21} \notin E(G)$ and hence $BA_2 \in E(G)$ (else color B by 2, $A_1$ by i, u by 1). Thus $b_{i2}$ has no 2-vertex. Then if $A_4a_{11} \notin E(G)$, color $b_{i2}$ by 2, $a_{12}$ by i, $A_4$ by 1, u by 4 and if $A_4a_{11} \in E(G)$, color $b_{i2}$ by 2, $A_4$ by i, u by 4, contradictions in both the cases. Hence $a_{21}a_{12} \notin E(G)$

$\Rightarrow a_{21}b_{i1} \in E(G)$ (else color $b_{i1}$ by 1, $a_{11}$ by i, $a_{21}$ by 1, $A_1$ by 2, u by 1).

Next let if possible $A_2a_{12} \in E(G)$. Then $b_{i1}$ is the only i-vertex of $A_2$ and $A_2B \notin E(G)$. Also $A_2$ is the only 2-vertex of $a_{12}$ and hence G has no other 2-vertex say $a_{22}$ (else $<a_{22}, a_{12}, A_1, a_{11}> = 4K_1$) $\Rightarrow Ba_{21} \in$

E(G) (else color B by 2, $A_1$ by i, u by 1). As $a_{21}$ has two 1-vertices and i-vertices, $a_{21}b_{i2} \notin E(G)$. Also as $A_2$ has two 1-vertices $A_2b_{i2} \notin E(G)$. Color $b_{i2}$ by 2, $a_{12}$ by i. If $a_{12}$ is the only 1-vertex of $A_4$, then color $A_4$ by 1, u by 4 and if $A_4a_{11} \in E(G)$, then color $a_{11}$ by 4, $A_4$ by 1, u by 4, contradictions in both the cases. Hence $A_2a_{12} \notin E(G)$. This proves **Claim 3**.

**Claim 4:** Whenever $A_1$ has a unique i-vertex B for $6 \leq i \leq \Delta-1$, either $A_2$ or $a_{21}$ has two i-vertices.
Let $b_{i1}$ be the only i-vertex $A_2$. Now $b_{i1}$ is not the only i-vertex of $a_{21}$ (else $<a_{21}, A_2, b_{i2}, B> = 4K_1$). Thus $a_{21}$ has two i-vertices. This proves **Claim 4.**

Now as $\Delta \geq 9$, and $A_1$ has at the most one repeat color, $A_1$ has at least two unique k-vertices for $k \in \{6, 7, ..., \Delta-1\}$. Let B, C be the unique i-vertex, k-vertex of $A_1$ resply for i, $k \in \{6, 7,..., \Delta-1\}$. Again as $a_{21}$ has two 1-vertices, each of $A_2$ and $a_{21}$ has at the most one other repeat color. By **Claim** 4, w.l.g. let $A_2$, $a_{21}$ have two i-vetices, k-vertices resply. $\Rightarrow A_2$, $a_{21}$ has a unique 4-vertex each. Similarly $A_4$, $a_{41}$ has a unique 2-vertex each. Now $A_2a_{41} \notin E(G)$ (else color $a_{41}$ by 2, $A_2$ by 4, u by 2). Also $A_2A_4 \notin E(G)$ (else color $A_4$ by 2, $A_2$ by 4, $a_{11}$ by 2, $a_{21}$ by 1, $A_1$ by 2, u by 1) $\Rightarrow a_{21}a_{41} \in E(G)$ (else $<a_{21}, a_{41}, A_2, A_4> = 4K_1$). As $a_{11}$ is the unique 1-vertex of $A_2$, color $a_{41}$ by 2, $a_{21}$ by 4, $a_{11}$ by 2, $A_2$ by 1, u by 2, a contradiction.

This proves **(C)**.

If $A_iA_j \notin E(G)$ ($1 \leq i, j \leq \Delta-2$), then as $\Delta-2 \geq 7$, by **(C)**, $\exists$ m ($1 \leq m \leq \Delta-2$) s.t. $A_iA_m$ $A_jA_m \in E(G)$. Also by **(B)**, $\exists$ maximum two such m's ($1 \leq m \leq \Delta-2$).

**Case 1.1:** $\exists$ i, j s.t. $A_iA_j \notin E(G)$ and $A_iA_k, A_iA_m, A_jA_k, A_jA_m \in E(G)$, $1 \leq i, j, k, m \leq \Delta-2$.
W.l.g. let i = 1, j = 2, k = 5, j = 6. Also by **(C)**, let $A_1A_4, A_2A_7 \in E(G)$. Then by **(B)**, $A_1A_7, A_2A_j \notin E(G)$. By **(A)**, w.l.g. let $A_1, A_2$ have two 5-vertices, 6-vertices resply. Clearly $A_4$ ($A_7$) is the unique 4-vertex (7-vertex) of $A_1$ ($A_2$). Also by **(C)**, $A_7$ is adjacent to at least one of $A_i$, $i \in \{3, 4, 6\}$ and if $A_7A_i \in E(G)$, $i \in \{3, 4, 6\}$, then $A_7$ has two i-vertices (else $A_7$, $A_1$ have a unique i-vertex $A_i$, a contradiction to **(A)**) and hence $A_2$ is the unique 2-vertex of $A_7$. Now $A_3A_1$ or $A_3A_2 \in E(G)$ (else by **(C)**, $A_3$ is adjacent to at least three of $A_4, A_5, A_6, A_7$ and either $A_3, A_1$ or $A_3, A_2$ have a common adjacent $A_i$ s.t. $A_i$ is their only i-vertex, a contradiction to **(A)**). W.l.g. let $A_3A_1 \in E(G)$. Again $A_3$ is the unique 3-vertex of $A_1$. Now $\exists$ 2-i paths from $A_2$ to $A_i$ (i = 1, 3, 4). Also as G is $4K_1$-free, $\exists$ at most two more 2-vertices $a_{21}, a_{22}$ and at least one of them say $a_{21}$ has two repeat colors from $\{1, 3, 4\}$.

**Case 1.1.1:** $a_{21}A_3, a_{21}A_4 \in E(G)$ and $a_{21}$ has two repeat colors 3, 4.
Then $a_{22}A_1 \in E(G)$ and $a_{22}$ has two 1-vertices and $a_{22}$ is the only 2-vertex of $A_1$. W.l.g. let $a_{22}$ have a unique 3–vertex (else $a_{22}$ has a color r missing in $N(a_{22})$. Color $a_{22}$ by r, $A_1$ by 2, u by 1). Then $a_{22}A_3 \notin E(G)$ (else color $a_{22}$ by 3, $A_3$ by 1, $A_1$ by 2, u by 3). Consider a 3-2 path T from $A_3$ to $A_2$ with $a_{31}$ being the 3-vertex of $A_2$ on T. As $a_{22}$ has a unique 3–vertex, clearly $a_{21}a_{31} \in E(G)$. Now $a_{22}a_{31} \in E(G)$ (else alter colors along $\{A_2, a_{31}, a_{21}, A_3\}$, color $A_1$ by 3, u by 1). Then G does not have a 3-vertex $a_{32} \notin \{A_3, a_{31}\}$ (else $<A_2, a_{22}, a_{32}, A_3> = 4K_1$). Now $A_7a_{31} \in E(G)$ (else $A_3$ is the only 3-vertex of both $A_1$ and $A_7$, contrary to **(A)**). But as $a_{31}$ has three 2-vertices, $A_7$ is its only 7-vertex. Also by **(C)**, $A_7$ is adjacent to at least one $A_j$ (j $\in \{3, 4, 6\}$ and has two j-vertices (else $A_j$ is the only j-vertex of $A_1$and $A_7$, contrary to **(A)**). Hence $A_2$ is the only 2-vertex of $A_7$. Then color $a_{31}$ by 7, $A_7$ by 2, $A_2$ by 3, u by 7, a contradiction. This proves **Case 1.1.1**.

**Case 1.1.2:** $A_3, A_4$ do not have a common adjacent 2-vertex.
W.l.g. let $a_{21}A_1, a_{21}A_3 \in E(G)$ and $a_{22}A_4 \in E(G)$. Then $a_{21}A_4, a_{22}A_3 \notin E(G)$. Clearly $a_{21}$ has two 1-vertices and 3–vertices and hence a unique 4-vertex. Let $a_{41}$ be the unique 4-vertex of $A_2$. Then as $\exists$ a 2-4 path S from $A_2$ to $A_4$, clearly $a_{22}a_{41} \in E(G)$. Now $a_{21}a_{41} \notin E(G)$ (else G does not have a 4-vertex $a_{42} \notin \{a_{41}, A_4\}$, as otherwise $<A_2, a_{21}, a_{42}, A_4> = 4K_1 \Rightarrow A_7a_{41} \in E(G)$ as otherwise $A_7$ and $A_1$ have a common unique 4-vertex $A_4$, a contradiction to **(A)**. But then color $A_7$ by 2, $A_2$ by 4, $a_{41}$ by 7, u by 7). Let $a_{42}$ be the unique 4-vertex of $a_{21}$. Then $a_{42}a_{22} \in E(G)$ (else alter colors along $\{A_4, a_{22}, a_{41}, A_2\}$, color $A_1$ by 4, u by 1). Thus $a_{22}$ has three 4-vertices and hence a unique i-vertex for $1 \leq i \leq \Delta-1$, i $\notin \{2, 4\}$. Now $A_4$ has a unique j-vertex for j = 1 or 3. Consider a 2-j path T from $A_2$ to $A_j$ and let $a_{j1}$ be the unique j-vertex of $A_2$. Clearly $a_{21}a_{j1} \in E(G)$. Again $a_{22}a_{j1} \notin E(G)$ (else G does not have a j-vertex $a_{j2} \notin \{a_{j1}, A_j\}$, as otherwise $<A_2, a_{22}, a_{j2}, A_j> = 4K_1 \Rightarrow A_7a_{j1} \in E(G)$. Color $A_7$ by 2, $A_2$ by j, $a_{j1}$ by 7, u by 7). Hence $\exists\ a_{j2}$ s.t. $a_{22}a_{j2} \in E(G)$ (else color $a_{22}$ by j, $A_4$ by 2, u by 4). Now clearly $a_{22}$ is the unique 2-vertex of $a_{j2}$ and vice versa $\Rightarrow A_4a_{j2} \in E(G)$ (else color $a_{22}$ by j, $a_{j2}$ by 2, $A_4$ by 2, u by 4). Clearly $A_4$ has two 1-vertices (else $A_1$ is its unique 1-vertex. Alter colors along $\{A_4, a_{22}, a_{41}, A_2\}$, color $A_1$ by 4, u by 1) $\Rightarrow$ j = 3 and $a_{32}$ is the unique 3-vertex of $A_4 \Rightarrow A_3A_4 \notin E(G)$. Then by **(C)**, $A_4$ is adjacent to at least two $A_k$s for k $\in$ \{5, 6, 7\}. Let $A_4A_m \in E(G)$, for m = 5 or 7. Then $A_m$ is the unique m-vertex of $A_4$ and $A_2$, a contradiction to **(A)**.

**Case 1.2:** $\forall$ i, j s.t. $A_iA_j \notin E(G)$, $A_i$, $A_j$ have only one common adjacent $A_k$ in N(u), $1 \leq i, j, k \leq \Delta-2$.
W.l.g. let i = 1, j = 2 and k = 3. By **(C)**, let $A_1A_m \in E(G)$ for m = 4, 5, $A_2A_l \in E(G)$ for l = 6, 7. **Let if possible $A_3A_4 \notin E(G)$**. Now $A_4$ is adjacent to at the most one of $A_6$, $A_7$ (else we get **Case 1.1** with $A_2$ and $A_4$) and hence by **(C)**, $A_4A_5 \in E(G)$. Also by **(C)**, w.l.g. let $A_4A_6 \in E(G)$. Again $A_3A_5$, $A_3A_6 \notin E(G)$ (else we get **Case 1.1** with $A_3$ and $A_4$) and hence by **(C)**, $A_3A_7 \in E(G) \Rightarrow A_5A_7 \notin E(G)$ (else we get **Case 1.1** with $A_3$ and $A_5$) and $A_5A_6 \in E(G)$. But then we get **Case 1.1** with $A_1$ and $A_6$, a contradiction. **Hence $A_3A_i \in E(G)$ for $4 \leq i \leq 7$**. Again $A_4A_5$, $A_6A_7 \in E(G)$ (else we get **Case 1.1** with $A_4$, $A_5$ or $A_6$, $A_7$). Also either all $A_i$ have two 3-vertices for $1 \leq i \neq 3 \leq 7$ or say $A_1$ has a unique 3-vertex. Again if $A_1$ has a unique 3-vertex, then $A_2$, $A_6$, $A_7$ all have two 3-vertices (else a contradiction to **(A)**). Hence w.l.g. let $A_1$, $A_4$, $A_5$ have two 3-vertices. As G is $4K_1$-free, G has at the most two 2-vertices say $a_{2i}$ (i = 1, 2). W.l.g. let $A_1a_{21} \in E(G)$. Now $a_{21}$ has at the most two repeat colors (else a color say r is missing in $N(a_{21})$. Color $a_{21}$ by r, $A_1$ by 2, u by 1). Also as $\exists$ i-2 paths from $A_i$ to $A_2$ for i = 1, 4, 5, either $a_{21}$ or $a_{22}$ has two j-vertices for j = 1, 4, 5. W.l.g. let $a_{21}$ have two repeat colors 1, 4 with $A_1a_{21}$, $A_4a_{21} \in E(G) \Rightarrow A_5a_{22} \in E(G)$ and $a_{22}$ has two 5-vertices. Again at least two of \{1, 4, 5\} are unique colors of $A_2$.

**Case 1.2.1.** $A_2$ has a unique 1-vertex and 5-vertex.
Let $A_2a_{11}$, $A_2a_{51} \in E(G)$. As $a_{21}$ has two repeat colors 1, 4, it has a unique 5-vertex and clearly as $\exists$ 2-5 path from $A_2$ to $A_5$, $a_{22}a_{51} \in E(G)$. Now $a_{21}a_{51} \notin E(G)$ (else G doesn't have a 5-vertex $a_{52} \notin \{A_5, a_{51}\}$ as otherwise $<A_5, a_{52}, A_2, a_{21}> = 4K_1$. As $a_{51}$ has three 2-vertices, $A_6$ is its only 6-vertex. Also $a_{51}$ is the only 5-vertex of $A_6$. Color $a_{51}$ by 6, $A_6$ by 5, u by 6) $\Rightarrow a_{21}a_{52} \in E(G)$. Also $a_{52}a_{22} \in E(G)$ (else color $A_2$ by 5, $a_{51}$ by 2, $a_{22}$ by 5, $A_5$ by 2, $A_1$ by 5, u by 1).. But then $a_{22}a_{11} \notin E(G)$ (else G doesn't have a 1-vertex $a_{12} \notin \{A_1, a_{11}\}$ as otherwise $<A_1, a_{12}, A_2, a_{22}> = 4K_1$ and $a_{11}$ has three repeat colors 2, 6, 7 with color say r missing in $N(a_{11})$. Color $a_{11}$ by r, $A_2$ by 1, u by 2) $\Rightarrow a_{11}a_{21} \in E(G)$. Let $a_{22}a_{12} \in E(G)$. Then $a_{22}$ ($a_{12}$) is the only 2-vertex (1-vertex) of $a_{12}$ ($a_{22}$). Color $a_{22}$ by 1, $a_{12}$ by 2, $A_5$ by 2, u by 5, a contradiction.

**Case 1.2.2.** $A_2$ has a unique 1-vertex and 4-vertex.
Let $A_2a_{11}$, $A_2a_{41} \in E(G)$. As $a_{22}$ has two 5-vertices, w.l.g. let $a_{22}$ have a unique 1-vertex. Then $a_{22}a_{11} \notin E(G)$ (else if $\exists\ a_{12}$, then $<A_1, a_{12}, A_2, a_{22}> = 4K_1$ and if $a_{12}$ doesn't exist, then $a_{11}$ has three repeat colors 2, 6, 7 and color say r is missing in $N(a_{11})$. Color $a_{11}$ by r, $A_2$ by 1, u by 2) $\Rightarrow a_{22}a_{12} \in E(G)$ and $a_{21}a_{11} \in E(G)$. Then $a_{22}$ ($a_{12}$) is the only 2-vertex (1-vertex ) of $a_{12}$ ($a_{22}$). Color $a_{22}$ by 1, $a_{12}$ by 2, $A_5$ by 2, u by 5, a contradiction.

This proves **Case 1**.

**Case 2:** In every $(\Delta-1)$-coloring of G-u, all vertices with unique colors in N(u) are adjacent.

Clearly $\Delta-1 \leq \omega$ and hence $\Delta-1 = \omega \geq 8 \Rightarrow <\bigcup_{i=1}^{\Delta-2} A_i>$ is a maximum clique in G-u and $\{X, Y\} = N(u) - \bigcup_{i=1}^{\Delta-2} A_i$.

**I.** At most two vertices in $\bigcup_{i=1}^{\Delta-2} A_i$ are non-adjacent to both X and Y.

Let if possible $A_1$, $A_2$, $A_3$ be non-adjacent to both X and Y. Then clearly ∃ a ($\Delta$-1)-vertex say Z in V(G) s.t. $ZA_i \in E(G)$ for i = 1, 2, 3. Moreover, as G is $4K_1$-free, Z is their only ($\Delta$-1)-vertex. If $A_i$ is the only i-vertex of Z for some i (1≤i≤3), then color $A_i$ by $\Delta$-1, Z by i, u by i, a contradiction. Hence Z has at least two i-vertices for i = 1, 2, 3. But then Z has some color r missing in N(Z). Color Z by r, $A_i$ by $\Delta$-1, u by i, a contradiction.

**II.** Every vertex $A_i$ of N(u) has at least one j-vertex j ≠ i (else color $A_i$ by j and u by i), 1≤i, j≤$\Delta$-2.

**III.** X (Y) has a k-vertex for every k = 1,.., $\Delta$-2.

Let if possible X not have a k-vertex. Also as $<u \cup \bigcup_{i=1}^{\Delta-2} A_i>$ is a maximum clique in G, ∃ i (1≤i≤$\Delta$-2) s.t. $YA_i \notin E(G)$. Then color X by k. Now i = k (else we get **Case 1** as Y and $A_i$ are unique vertices in N(u)). As $\Delta \geq 9$ and each of Y and $A_i$ has at the most one repeat color, clearly ∃ j (1≤j≤$\Delta$-2) s.t. $A_j$ is the only j-vertex of both Y and $A_i$. Also $A_j$ has either a unique i-vertex $A_i$ or ($\Delta$-1)-vertex Y. Color Y and $A_i$ by j, $A_j$ by i ($\Delta$-1), u by $\Delta$-1 (j), a contradiction.

**IV.** X (Y) is adjacent to at least ω-5 vertices in $\bigcup_{i=1}^{\Delta-2} A_i$.

Let if possible X be non-adjacent to $A_i$, i = 1,.., 5. By **I**, w.l.g. let $YA_i \in E(G)$ for i = 1, 2, 3. Also let $YA_k \notin E(G)$ for some k ≥ 4. By **II** and **III**, Y and $A_k$ each has at the most one repeat color and hence w.l.g. let $A_1$ be the unique 1-vertex of Y and $A_k$. Now $A_1$ has two ($\Delta$-1)-vertices (else color Y and $A_k$ by 1, $A_1$ by $\Delta$-1, u by k) ⇒ $A_k$ is the unique k-vertex of $A_1$. Then color Y and $A_k$ by 1, $A_1$ by k and we get **Case 1** with two non-adjacent, unique vertices X, $A_1$, a contradiction.

**V.** X (Y) is not the only ($\Delta$-1)-vertex of any $A_i$.

Let if possible X be the only ($\Delta$-1)-vertex of some $A_i$. By **IV**, ∃ k, j s.t. $XA_k$, $XA_j \in E(G)$. Also let $XA_m \notin E(G)$ for some m. If $A_i$ is the only i-vertex of X and $A_m$, then color X, $A_m$ by i, $A_i$ by $\Delta$-1, u by m, a contradiction. Hence let $A_i$ be not the only i-vertex of either X or $A_m$. As X and $A_m$ have at the most one repeat color, w.l.g. let $A_k$ be the only k-vertex of X and $A_m$. Again if X is the only ($\Delta$-1)-vertex of $A_k$, then as before we get a contradiction. Hence let $A_k$ have two ($\Delta$-1)-vertices. But then color $A_k$ by i, $A_i$ by $\Delta$-1, X by k, $A_m$ by k, u by m, a contradiction.

By **IV**, w.l.g. let $XA_k \in E(G)$ for k = 1, 2, 3 and $XA_4 \notin E(G)$. Also w.l.g. let $A_1$ be the only 1-vertex of X and $A_4$. By **V**, $A_1$ has two ($\Delta$-1)-vertices. If any $A_i$ (1≤i≤$\Delta$-2, i≠4) is non-adjacent to Y, then as before by coloring X, $A_4$ by 1 and $A_1$ by 4, we get **Case 1** and hence $YA_k \in E(G)$ for every k ≠ 4. Similarly $XA_k \in E(G)$ for every k ≠ 4. As $\Delta \geq 9$, ∃ i s.t. $A_i$ is the only i-vertex of X, Y and $A_4$. Color X, Y, $A_4$ by i, $A_i$ by $\Delta$-1, u by 4, a contradiction.

This proves **Case 2** and completes the proof of the Main Result.

**References**

[1] "Linear Chromatic Bounds for a Subfamily of 3K1-free Graphs",S. A. Choudum, T. Karthick, M. A. Shalu,Graphs and Combinatorics 24:413–428, 2008

[2] "On the divisibility of graphs", Chinh T. Hoang, Colin McDiarmid, Discrete Mathematics 242, 145–156, 2002

[3] "ω, $\Delta$, and χ", B.A. Reed, J. Graph Theory 27, pp. 177-212, 1998

[4] "Some results on Reed's Conjecture about ω, $\Delta$ and χ with respect to α", Anja Kohl, Ingo Schiermeyer, Discrete Mathematics 310, pp. 1429-1438, 2010

[5] `"On colouring the nodes of a network", Proc. Cambridge Philosophical Society, Math. Phys. Sci., 37 (1941), 194–197

[6] O. V. Borodin and A. V. Kostochka, "On an upper bound of a graph's chromatic number, depending on the graph's degree and density", JCTB 23 (1977), 247--250.

[7] B. A. Reed, "A strengthening of Brooks' Theorem", J. Comb. Theory Ser. B, 76 (1999), 136–149.

[8] D. W. Cranston and L. Rabern, "Coloring claw-free graphs with $\Delta$-1 colors" *SIAM J. Discrete Math.*, 27(1) (1999), 534–549.

[9] J.A. Bondy and U.S.R. Murty. Graph Theory, volume 244 of Graduate Text in Mathematics. Springer, 2008.